\documentclass[8pt]{article}

\usepackage{amsmath}

\usepackage{amssymb}

\usepackage{epsfig}
\numberwithin{equation}{section}
\newtheorem{thm}{Theorem}[section]
\newtheorem{lem}{Lemma}[section]

\newtheorem{defi}{Definition}[section]
\newtheorem{rem}{Remark}[section]
\newtheorem{asm}{Assumption}[section]

\begin{document}
\title{Constrained LQ problem with a random jump and application to portfolio selection \footnotemark[1]}
\author{Yuchao Dong\footnotemark[2]}
\maketitle
\renewcommand{\thefootnote}{\fnsymbol{footnote}}

\footnotetext[1]{This research was supported by the National Natural Science Foundation of China (Grants \#10325101 and \#11171076), and by Science and Technology Commission, Shanghai Municipality (Grant No.14XD1400400). The author would like to thank his advisor, Prof. Shanjian Tang from Fudan University, for helpful comments and discussions.}

\footnotetext[2]{Department of Finance and Control Sciences, School of Mathematical Sciences, Fudan University, Shanghai 200433, China. yuchaodong13@fudan.edu.cn (Yuchao Dong).}
\begin{abstract}
 In this paper, we consider a constrained stochastic linear-quadratic (LQ) optimal control problem where the control is constrained in a closed cone. The state process is  governed by a controlled SDE with random coefficients. Moreover, there is a random jump of the state process. In mathematical finance, the random jump often represents the default of a counter party. Thanks to the It\^o-Tanaka formula, optimal control and optimal value can be obtained by solutions of a system of backward stochastic differential equations (BSDEs). The solvability of the BSDEs is obtained by solving a recursive system of BSDEs driven by the Brownian motions. We also apply the result to the mean variance portfolio selection problem in which the stock price can be affected by the default of a counterparty.
\end{abstract}
\begin{section}{Introduction}
~~~~Linear-quadratic (LQ) problem is an important optimal control problem. The feature of such a problem is that the dynamic of the system is linear in the state and control variables and the cost functional is quadratic in both of them. It was first considered by Kalman \cite{Kal1} (for the deterministic control of ordinary differential equations, i.e. ODEs) and then extended to various situations, for example stochastic LQ problems. One important application of stochastic LQ optimal control theory is the continuous-time version of Markowitz's mean-variance portfolio selection problem, which is one fundamental problem in the mathematical finance.

It is well-known that one can give in explicit forms the optimal state feedback control and the optimal value via the celebrated Riccati equation. In the deterministic case or the stochastic case with deterministic coefficients, the Riccati equation is an ODE in the space of symmetric matrixes. When the coefficients are random, the Riccati equation becomes a backward stochastic differential equation (BSDE). The theory of BSDEs was pioneered by Pardoux and Peng \cite{PP1}. It is closely related to the optimal control theory. See Yong and Zhou \cite{YZ1} on this subject. For Riccati equations, the solvability is a very hard problem. Under some standard assumptions of the coefficients, it is solved by Tang \cite{Tang1,Tang2} by two different approaches. For more details on this subject, see \cite{Du_2015,HZ1,Tang1,Tang2} and the references therein.

In this paper, we consider the stochastic LQ problems with a random jump. Note that similar problems have also been considered by \cite{HO1}, \cite{Meng1} and \cite{MT1}. Our problem is different from theirs from two apsects. One is that our system only has at most one jump. In mathematical finance, this random jump represents the default, so sometimes we just call it the default time. In a financial market, we know that the default of one firm has usually important influences on the others. This has been shown clearly in the financial crisis. While the controlled processes considered in those papers mentioned above are driven by a Poisson random measure, their systems can have even infinitely many jumps. The other difference is that the control in our problem is constrained in a closed cone. In the mean-variance problem, this means that there are some restrictions on the trading strategy of the investor. In this paper, we shall consider the mean-variance portfolio selection problem for an investor who invests in a risky asset exposed to a counterparty risk. The investor is also not allowed to short sell. Thus we have to solve a constraint LQ problem with a random jump. We only consider the problem for the case that the state variable is scalar-valued. How to solve it in the multi-dimensional case is still a problem, but the scalar-valued case is sufficient to cover many important practical applications especially in the financial area.

To get the optimal control and the optimal value, we must first get the Riccati equation. Note that, due to the constraint, the value function is no longer quadratic with respect to the initial value.  But one can easily show that the value function $V$ is positive homogeneous since the control is constrained in a closed cone. That is
\[V(t,x)=\frac{1}{2}P_tx^{+,2}+\frac{1}{2}N_tx^{-,2},\]
where $P$ and $N$ satisfies the following BSDEs:
\begin{equation}
\begin{split}
d{P_t} =&  - \{ 2({A_t}{P_{t - }} - {{\tilde \lambda }_t}{E_t}{P_{t - }}) + \left\langle {{C_t},{Z_t}} \right\rangle  + {Q_t} \\
&+ {h^ + }(t,{P_{t - }},{Z_t},{{\bar Z}_t},{\bar \Lambda _t} + {N_{t - }})\} dt + {Z_t}d{W_t} + {{\bar Z}_t}d{M_t}\\
{P_T} =& G,
\end{split}
\end{equation}
\begin{equation}
\begin{split}
d{N_t} =&  - \{ 2({A_t}{N_{t - }} - {{\tilde \lambda }_t}{E_t}{N_{t - }}) + \left\langle {{C_t},{Z_t}} \right\rangle  + {Q_t} \\
&+ {h^ - }(t,{N_{t - }},{\Lambda _t},{{\bar \Lambda }_t},{{\bar Z}_t} + {P_{t - }})\} dt + {\Lambda _t}d{W_t} + {{\bar \Lambda }_t}d{M_t}\\
{N_T} =& G,
\end{split}
\end{equation}
Thus we are still able to get a system of BSDEs, sometimes called extended Riccati equation, that characterizes the optimal control and the optimal value. We can see that the BSDEs are coupled and have a random jump. Note that multi-dimensional backward Riccati equations have also been considered by  K. Mitsui and Y. Tabata \cite{MT1}. But their equations are multi-dimensional because the state processes in \cite{MT1} are multi-dimensional. To solve such equations, we use the method originated by Ankirchner et al \cite{ABE1} and further developed by Kharroubi and Lim \cite{KL1}. Through the decomposition of processes with respect to the progressive enlargement of filtrations, we link the BSDEs we want to solve with a family of Brownian BSDEs. By proving the solvability of the Brownian BSDEs, we are able to solve the original BSDEs. If there is no jump, the equations will be decoupled and this is the exact equation considered by Hu and Zhou \cite{HZ1}.

The rest of the paper is organized as follows. In Section 2, we formulate the problem. In Section 3, we derive the form of the extend Riccati equations and prove its solvability in two cases. In Section 4, we give the state feedback optimal control and the optimal value via the Riccati equations. The application to mean-variance problem is in Section 5.
\end{section}
\begin{section}{The model and Assumptions}
~~~~In this paper, we assume throughout that $(\Omega , \mathcal {F}, \mathbb {P})$ is a given probability space and that
$W_t$ is a $k$-dimensional standard Brownian motion defined on this space with $W_0=0$. Let $\{\mathcal {F}_t\}$ be the augmentation of $\sigma \{W_s|0 \le s \le t\}$. In addition, let $\tau$ be a random time. Define
\[\mathcal {G}_t=\mathop  \cap \limits_{s \ge t} \mathcal {F}_s \vee \sigma (1_{\tau \le s}),\]which is the smallest filtration containing $\{\mathcal {F}_t\}$ that makes $\tau$ a stopping time and satisfies the usual condition.

Throughout this paper, we denote the inner product of $\mathcal {R}^m$ by $\left\langle\cdot,\cdot\right\rangle$. If $M \in \mathcal {S}^n$ is positive (positive semi-) definite, we write $M>$($\ge$)$0$. Let $\mathbb {F}=\{\mathcal {F}_t\}$ be a filtration. Denote by $\mathcal P(\mathbb {F})$ the $\sigma$-field of $\mathbb {F}$-predictable measurable subsets of $\Omega \times \mathcal {R}_+$. Suppose that $f$ is a $\mathcal {R}^n$-valued square integrable process (i.e. $E[\int_0^T\left|f_s\right|^{2}ds]< \infty$). If $f$ is  $\mathbb F$-adapted, we shall write $f\in \mathcal {L}_{\mathbb F}^2([0,T],\mathcal R^n)$; if $f$ is  $\mathbb F$-predictable, then $f \in \mathcal {L}^2([0,T]\times \Omega,\mathcal {P}(\mathbb F),\mathcal R^n)$. Similarly, denote by $\mathcal {L}_{\mathbb F}^{\infty}([0,T],\mathcal R^n)(\mathcal {L}^{\infty}([0,T]\times \Omega,\mathcal {P}(\mathbb F),\mathcal R^n)$) the set of all bounded adapted (predictable) processes. Furthermore denote by $\mathcal {S}_{\mathbb F}^{\infty}([0,T],\mathcal R^n)$ ($(\mathcal {S}_{\mathbb F}^{2}([0,T],\mathcal R^n)$) the set of processes that belong to $\mathcal {L}_{\mathbb F}^{\infty}([0,T],\mathcal R^n)$ ($\mathcal {L}_{\mathbb F}^{2}([0,T],\mathcal R^n)$) with continuous paths. These definitions generalize in the obvious way to the case when $f$ is $\mathcal {R}^{n \times m}$- or $\mathcal {S}^n$-valued. Moreover, we say that $N \in \mathcal {L}^2([0,T]\times \Omega,\mathcal {P}(\mathbb F),\mathcal S^n)$ is positive (positive semi-) definite, which is sometimes denote by $N >(\ge)0$, if $N(t,\omega)>(\ge)0$ for a.e. $t \in [0,T]$ and $\mathbb P$-a.s., and say that $N$ is uniformly positive definite if $N(t,\omega) \ge cI_n$ for a.e. $t \in [0,T]$ and $\mathbb P$-a.s. with some deterministic constant $c$, where $I_n$ is the $n$-dimensional identity matrix. Finally, for any real number we define $x^+=\max\{x,0\},x^-=\max\{-x,0\},x^{+,2}=(x^+)^2$ and $x^{-,2}=(x^-)^2$.

In the sequel, we shall make the following assumptions on the random time $\tau$.  For any $t \in [0,T]$, the conditional distribution of $\tau$ under $\mathcal {F}_t$ admits a density with respect to Lebsegue measure, i.e. there exists an $\mathcal {F}_t \otimes \mathcal {B}(\mathcal R_+)$-measurable positive function $(\omega,\theta) \to \alpha_t(\theta)$ such that
\begin{equation}
P[\tau \in d\theta|\mathcal {F}_t]=\alpha_t(\theta)d\theta.
\end{equation}
Note that for any $\theta \ge 0$, the process $\{\alpha_t(\theta),0\le t\le T\}$ is a $\mathbb F$-martingale. Moreover we assume that the family of densities satisfy $\alpha _T(t)=\alpha_t(t)$ for all $0 \le t \le T$.
\begin{rem}In the finance, the random time $\tau$ usually represents the default of a counterparty. The density hypothesis is usually used in the theory of enlargement of filtrations. It was  introduced in the notes of  Jeulin and Yor \cite{JY1} and recently adopted by El Karoui et al \cite{EJJ1} for credit risk modelling. Note that we have $P[\tau >t||\mathcal{F}_t]=P[\tau >t||\mathcal{F}_T]$. This is related to the so-called immersion hypothesis meaning that any square integrable $\mathbb F$-martingale is a square integrable $\mathbb G$-martingale.\end{rem}
Let $L_t=1_{\{\tau \le t\}}$. Then $L$ is a $\{\mathcal {G}_t\}$-submartingale. We will have the following assumption.
\begin{asm}
There exists an $\mathbb {F}$-predicable bounded nonnegative process $\lambda$ such that
\begin{equation}
M_t:=L_t-\int _0^t(1-L_{s-})\lambda _sds
\end{equation}
is a martingale with respect to $\{\mathcal {G}_t\}$.
\label{asm-2.1}\end{asm}
We also define $\tilde \lambda _s:=(1-L_{s-})\lambda _s$.  Let us mention that $\lambda$ can be explicitly expressed by the conditional density (see \cite{EJJ1}). In fact, $\lambda _t=\alpha_t(t)/G(t)$, where $G(t)=P[\tau >t||\mathcal{F}_t]=\int_t^{\infty}\alpha_t(\theta)d\theta.$
Now we give one example that the assumption holds.\\
\begin{bfseries} Example\end{bfseries}  Let $\beta$ be a bounded nonnegative $\{\mathcal {F}_t\}$-predictable process such that
\[\int_0^{\infty}\beta_sds=+\infty \text{\quad a.s.} \]
and $\Theta$ an exponential distributed random variable that is independent of the Brownian motion $W$. Define the random time
\[\tau = \inf\{t;\int _0^t \lambda _s ds > \Theta\}.\]
Then one can show that $\alpha_t(\theta)=E[\beta_\theta e^{-\int_0^\theta\beta_sds}|\mathcal F_t]$ and $\lambda_t=\beta_t$. Thus the assumption is satisfied. We refer the readers to the monograph of Jeanblanc et al \cite{JMC1} for the details and its application in mathematical finance.
\begin{rem}
Let $\phi$ be a $\{\mathcal {G}_t\}$-predictable process. Then it can be represented as
\begin{equation}
\phi_t=\phi ^{0}1_{\{t \le \tau\}}+ \phi _t^1(\tau)1_{\{\tau < t\}},
\end{equation}
where $\phi ^0$ is $\mathbb F$-predictable and $\phi^1$ is $\mathcal {P}(\mathbb {F})\otimes \mathcal {B}(\mathcal R)$-measurable.
\label{rm-2.1}\end{rem}
Consider the following controlled linear SDE:
\begin{equation}
\begin{split}
d{X_s}& = \{ {A_s}{X_{s - }} + {B_s}{u_{s}}\} ds + \{ {C_s}{X_{s - }} + {D_s}{u_{s}}\} d{W_s}\\
         &\quad + \{ {E_s}{X_{s - }} + {F_s}{u_{s}}\} d{M_s},\qquad t \le t \le T,\\
{X_t}& = x.
\end{split}
\label{csde}\end{equation}
The coeifficients $A.B,C,D,E,F$ are $\{\mathcal {G}_t\}$-predictable processes, and $x \in \mathcal{R}$ is a nonrandom scalar. Precise assumptions on these coefficients will be specified below. Let $\Gamma \subseteq \mathcal {R}^m$ be a given closed cone. A typical example of such a cone is $\Gamma=\mathcal {R}_+^m$. The class of  admissible controls is the set
$\mathcal {U}:=\mathcal {L}^2([0,T]\times \Omega, \mathcal {P}(\mathbb G),\Gamma),$ i.e. the square integrable $\Gamma$-valued $\{\mathcal{G}_t\}$-predictable processes. The cost is given by
\begin{equation}
J(t,x,u) := E^{\mathcal {G}_t}[\frac{1}{2}GX_T^2 + \frac{1}{2}\int_t^T {{Q_r}} X_r^2 + \left\langle {{R_r}{u_r},{u_r}} \right\rangle dr]
\end{equation}
The optimal control problem is to minimize the cost functional over all admissible controls. Define the value function by
\[V(t,x)=\mathop {\text{essinf}}\limits_{u \in \mathcal U}J(t,x,u).\]
We have the following assumptions on the coefficients.
\begin{asm}
\[A,Q \in \mathcal {L}^{\infty }([0,T] \times \Omega ,\mathcal {P}(\mathbb G),\mathcal R),B \in \mathcal {L}^{\infty }([0,T] \times \Omega ,\mathcal {P}(\mathbb G),\mathcal R^m),\]
\[C \in \mathcal {L}^{\infty }([0,T] \times \Omega ,\mathcal {P}(\mathbb G),\mathcal R^k),D \in \mathcal {L}^{\infty }([0,T] \times \Omega ,\mathcal {P}(\mathbb G),\mathcal R^{k \times m}),\]
\[E \in \mathcal {L}^{\infty }([0,T] \times \Omega ,\mathcal {P}(\mathbb G),\mathcal R), F \in \mathcal {L}^{\infty }([0,T] \times \Omega ,\mathcal {P}(\mathbb G),\mathcal R^m),\]
\[R \in \mathcal {L}^{\infty }([0,T] \times \Omega ,\mathcal {P}(\mathbb G),\mathcal S^m),G \in \mathcal {L}^{\infty}( \mathcal G_T,\mathcal R).\]
\[E(t,\omega) \ge -1,dtd\mathbb{P}\text{-a.s.}.\]
\label{asm-2.2}\end{asm}
By the Remark \ref{rm-2.1}, we will have the following decompositions of the processes
\begin{equation}
i(t)=i^0(t)1_{\{t \le \tau\}}+i^1_t(\tau)1_{\{\tau <t \}},
\end{equation}
where $i ^0$ is $\mathbb F$-predictable and $i^1$ is $\mathcal {P}(\mathbb {F})\otimes \mathcal {B}(\mathcal R)$-measurable for $i$=A,B,C,D,E,F,R,Q. And
\begin{equation}
G=G^01_{\{t \le \tau\}}+G^1(\tau)1_{\{\tau <t \}},
\end{equation}
where $G^0$ is $\mathcal {F}_T$-measurable and $G^1$ is $\mathcal  {F}_T \otimes \mathcal {B}(\mathcal R)$-measurable.

\end{section}

\begin{section}{Existence of solutions for the stochastic Riccati Equations}
\subsection{the Form of the Riccati Equations}
~~~~In this section, we will prove the existence of solutions for the extended stochastic Riccati Equations. First of all, we shall derive the formation of the Riccati equations. Note that the admissible controls are $\Gamma$-valued and $\Gamma$ is a closed cone. It means that for any $u \in \mathcal{U}$ and $c \ge 0$, $cu$ also belongs to $\mathcal {U}$. Since the controlled SDE is linear and the cost functional is quadratic, it is obvious that the value function $V$ is positive homogeneous, i.e. $V(t,cx)=c^2V(t,x)$ for all $c \ge 0$. Hence $V$ is of the following form
\begin{equation}
V(t,x)=\frac{1}{2}P_tx^{+,2}+\frac{1}{2}N_tx^{-,2}.
\end{equation}
Assume that both $P$ and $N$ are semimartingales with the following decompositions
\begin{align}
dP_t=f_tdt+Z_tdW_t+\bar Z_tdM_t,P_T=G,\label{pbsde}\\
dN_t=g_tdt+\Lambda_tdW_t+\bar \Lambda_tdM_t,N_T=G.\label{nbsde}
\end{align}
Given any $u \in \mathcal {U}$, $X$ is the associated solution of (\ref{csde}). By the It\^o-Tanaka formula, we have
\begin{equation*}
\begin{split}
\frac{1}{2}X_s^{ + ,2} =& \frac{1}{2}X_t^{ + ,2} + \int_t^s {X_r^ + } d{X_r} + \int_t^s {\frac {1}{2}{1_{\{ {X_{r - }} \ge 0\} }}} d{\left\langle {{X^c}} \right\rangle _r} \\
&+ \sum\limits_{t < r \le s} {(\frac{1}{2}X_r^{ + ,2} - } \frac{1}{2}X_{r - }^{ + ,2} - X_{r - }^ + \Delta {X_r}).
\end{split}
\end{equation*}
Note that $X$ only has a jump at the time $\tau$, i.e.
\[\Delta {X_s} = \left\{ {\begin{array}{*{20}{c}}
0\\
{{E_\tau }{X_{\tau  - }} + {F_\tau }{u_\tau }}
\end{array}\begin{array}{*{20}{c}}
{\text{otherwise}}\\
{s = \tau }
\end{array}} \right..\]
Hence we get that
\[\sum\limits_{t < r \le s} {(\frac{1}{2}X_r^{ + ,2} - } \frac{1}{2}X_{r - }^{ + ,2} - X_{r - }^ + \Delta {X_r}) = \int_t^s {f_r^ + ({X_{r - }},{F_r}{u_r})}  - X_{r - }^ + ({E_r}{X_{r - }} + {F_r}{u_r})d{L_r},\]
where $f_t^ + (x,y) = {\frac {1}{2}(x + {E_t}x + y)^{ + ,2}} - \frac{1}{2}{x^{ + ,2}}.$\\
Thus
\begin{equation*}
\begin{split}
&\frac{1}{2}X_s^{ + ,2}\\
 = &\frac{1}{2}X_t^{ + ,2} + \int_t^s {} \{ {A_r}X_{r - }^{ + ,2} + X_{r - }^ + {B_r}{u_r}\\
  &+ {{\tilde \lambda }_r}(f_r^ + ({X_{r - }},{F_r}{u_r}) - X_{r - }^ + ({E_r}{X_{r - }} + {F_r}{u_r}))+ \frac {1}{2}{1_{\{ {X_{r - }} \ge 0\} }}{\left| {{C_r}{X_{r - }} + {D_r}{u_r}} \right|^2}\} dt \\
  &+ \int_t^s {{C_r}X_{r - }^{ + ,2} + X_{r - }^ + {D_r}{u_r}} d{W_r} + \int_t^s {f_r^ + ({X_{r - }},{F_r}{u_r})d{M_r}}.
\end{split}
\end{equation*}
By (\ref{pbsde}) and It\^o formula again,
\begin{equation}
\begin{aligned}
\frac{1}{2}{P_s}X_s^{ + ,2}= &\frac{1}{2}{P_t}X_t^{ + ,2} + m_s \\
&+\int_t^s {} \bigg\{ {1_{\{ {X_{r - }} \ge 0\} }}\big\{ \frac{1}{2}{f_r}X_{r - }^2 + {P_{r - }}({A_r}X_{r - }^{ + ,2} + X_{r - }^ + {B_r}{u_r} \\
&- {{\tilde \lambda }_r}(X_{r - }^ + ({E_r}{X_{r - }} + {F_r}{u_r})) +\frac{1}{2} {\left| {{C_r}{X_{r - }} + {D_r}{u_r}} \right|^2}) \\
&+ \left\langle {{C_r}X_{r - }^2 + {X_{r - }}{D_r}{u_r},Z_r} \right\rangle \big\}  + ({{\bar Z}_r} + {P_{r - }}){{\tilde \lambda }_r}f_r^ + ({X_{r - }},{F_r}{u_r})\bigg\} dr,
\end{aligned}\label{ppdec}
\end{equation}
where $m_s$ is the local martingale part:
\[
\begin{split}
{m_s} = &\int_t^s {\big \{({C_r}X_{r - }^{ + ,2} + X_{r - }^ + {D_r}{u_r}){P_{r - }} + \frac{1}{2}X_{r - }^{ + ,2}{Z_r}\big \}d{W_r}} \\
 &+ \int_t^s {\big \{({P_{r - }} + {{\bar Z}_r})f_r^ + ({X_{r - }},{F_r}{u_r}) + \frac{1}{2}X_{r - }^{ + ,2}{{\bar Z}_r}\big \}d{M_r}}.
\end{split}
\]
Similarly, we also have
\begin{equation}
\begin{aligned}
\frac{1}{2}{N_s}X_s^{ - ,2} =&\frac{1}{2}{N_t}X_t^{ -,2} + n_s \\
&+\int_t^s {} \bigg\{ {1_{\{ {X_{r - }} \le 0\} }}\big\{ \frac{1}{2}{g_r}X_{r - }^2 + {N_{r - }}({A_r}X_{r - }^{ - ,2} - X_{r - }^ - {B_r}{u_r} \\
&+{{\tilde \lambda }_r}({X_{r - }^-}({E_r}{X_{r - }} - {F_r}{u_r}))+ \frac{1}{2}{\left| {{C_r}{X_{r - }} + {D_r}{u_r}} \right|^2}) \\
&+ \left\langle {{C_r}X_{r - }^2 - {X_{r - }^-}{D_r}{u_r},{\Lambda _r}} \right\rangle \big\}  + ({\bar \Lambda _r} + {N_{r - }}){{\tilde \lambda }_r}f_r^ - ({X_{r - }},{F_r}{u_r})\bigg\} dr,
\end{aligned}\label{npdec}
\end{equation}
where $f_t^ - (x,y) = {\frac{1}{2}(x + {E_t}x + y)^{ - ,2}} - \frac{1}{2}{x^{ - ,2}}$ and
\begin{equation*}
\begin{split}
{n_s} =& \int_t^s {({C_r}X_{r - }^{ - ,2} - X_{r - }^ - {D_r}{u_r}){N_{r - }} + \frac{1}{2}X_{r - }^{ - ,2}{\Lambda _r}d{W_r}} \\
 &+ \int_t^s {({N_{r - }} + {{\bar \Lambda }_r})f_r^ - ({X_{r - }},{F_r}{u_r}) + \frac{1}{2}X_{r - }^{ - ,2}{{\bar \Lambda }_r}d{M_r}}
\end{split}
\end{equation*}
Combining (\ref{ppdec}) and (\ref{npdec}) and letting $s=T$,
\begin{equation}
\begin{split}
&\frac{1}{2}GX_T^2 + \frac{1}{2}\int_t^T {{Q_r}} X_r^2 + \left\langle {{R_r}{u_r},{u_r}} \right\rangle dr\\
= &V(t,x) + m_T+n_T  \\
&+\int_t^T {} \{ {1_{\{ {X_{r - }} \ge 0\} }}\{ \frac{1}{2}{f_r}X_{r - }^2 + {P_{r - }}({A_r}X_{r - }^{ + ,2} + X_{r - }^ + {B_r}{u_r} - {{\tilde \lambda }_r}({X_{r - }}({E_r}{X_{r - }} + {F_r}{u_r}))\\
 &+ \frac{1}{2}{\left| {{C_r}{X_{r - }} + {D_r}{u_r}} \right|^2}) + \left\langle {{C_r}X_{r - }^2 + {X_{r - }}{D_r}{u_r},{\Lambda _r}} \right\rangle \}  + ({{\bar Z}_r} + {P_{r - }}){{\tilde \lambda }_r}f_r^ + ({X_{r - }},{F_r}{u_r})\\
 &+ ({\bar \Lambda _r} + {N_{r - }}){{\tilde \lambda }_r}f_r^ - ({X_{r - }},{F_r}{u_r}) + \frac{1}{2}{Q_r}X_r^2 + \frac{1}{2}\left\langle {{R_r}{u_r},{u_r}} \right\rangle \} dr\\
 &+ \int_t^T {} \{ {1_{\{ {X_{r - }} \le 0\} }}\{ \frac{1}{2}{g_r}X_{r - }^2 + {N_{r - }}({A_r}X_{r - }^{ - ,2} - X_{r - }^ - {B_r}{u_r} + {{\tilde \lambda }_r}({X_{r - }^-}({E_r}{X_{r - }} - {F_r}{u_r}))\\
 &+ \frac{1}{2}{\left| {{C_r}{X_{r - }} + {D_r}{u_r}} \right|^2}) + \left\langle {{C_r}X_{r - }^2 - {X_{r - }^-}{D_r}{u_r},{\Lambda _r}} \right\rangle \}  + ({{\bar Z}_r} + {P_{r - }}){{\tilde \lambda }_r}f_r^ + ({X_{r - }},{F_r}{u_r})\\
 &+ ({\bar \Lambda _r} + {N_{r - }}){{\tilde \lambda }_r}f_r^ - ({X_{r - }},{F_r}{u_r})+\frac{1}{2}{Q_r}X_r^2 + \frac{1}{2}\left\langle {{R_r}{u_r},{u_r}} \right\rangle \} dr.
\end{split}\label{value.fun}
\end{equation}
We denote that
\begin{align*}
g^+(r,X_r,u_r):=&{P_{r - }}({A_r}X_{r - }^{ + ,2} + X_{r - }^ + {B_r}{u_r} - {{\tilde \lambda }_r}({X_{r - }}({E_r}{X_{r - }} + {F_r}{u_r}))\\
 &+ \frac{1}{2}{\left| {{C_r}{X_{r - }} + {D_r}{u_r}} \right|^2}) + \left\langle {{C_r}X_{r - }^2 + {X_{r - }}{D_r}{u_r},{\Lambda _r}} \right\rangle \}  \\
 &+ ({{\bar Z}_r} + {P_{r - }}){{\tilde \lambda }_r}f_r^ + ({X_{r - }},{F_r}{u_r})\\
 &+ ({\bar \Lambda _r} + {N_{r - }}){{\tilde \lambda }_r}f_r^ - ({X_{r - }},{F_r}{u_r}) + \frac{1}{2}{Q_r}X_r^2 + \frac{1}{2}\left\langle {{R_r}{u_r},{u_r}} \right\rangle
\end{align*}
and
\begin{align*}
g^-(r,X_r,u_r):=&{N_{r - }}({A_r}X_{r - }^{ - ,2} - X_{r - }^ - {B_r}{u_r} + {{\tilde \lambda }_r}({X_{r - }^-}({E_r}{X_{r - }} - {F_r}{u_r}))\\
 &+ \frac{1}{2}{\left| {{C_r}{X_{r - }} + {D_r}{u_r}} \right|^2}) + \left\langle {{C_r}X_{r - }^2 - {X_{r - }^-}{D_r}{u_r},{\Lambda _r}} \right\rangle \}  \\
 &+ ({{\bar Z}_r} + {P_{r - }}){{\tilde \lambda }_r}f_r^ + ({X_{r - }},{F_r}{u_r})\\
 &+ ({\bar \Lambda _r} + {N_{r - }}){{\tilde \lambda }_r}f_r^ - ({X_{r - }},{F_r}{u_r})+\frac{1}{2}{Q_r}X_r^2 + \frac{1}{2}\left\langle {{R_r}{u_r},{u_r}} \right\rangle.
\end{align*}
Since $V$ is the value function, the integrand should always be positive. For some admissible control $u$, if the integrand is zero and the local martingale part is in fact a martingale, then taking conditional expectation, we have that it will be the optimal control. Hence we must have that
\[f_rX_r^{+,2} \ge -g^+(r,X_r,u_r)=-g^+(r,1,\frac {u_r}{X_r^+})X_r^{+,2}.\]
Noting that $\Gamma$ is a close cone, we have $\frac {u_r}{X_r^+} \in \Gamma$, thus $f_t$ should satisfy
\[f_t=-\inf \limits_{v \in \Gamma}g(t,1,v).\]
With a similar discussion, we see that $P$ and $N$ should be the solutions of the following system of BSDEs:
\begin{equation}
\begin{split}
d{P_t} =&  - \{ 2({A_t}{P_{t - }} - {{\tilde \lambda }_t}{E_t}{P_{t - }}) + \left\langle {{C_t},{Z_t}} \right\rangle  + {Q_t} \\
&+ {h^ + }(t,{P_{t - }},{Z_t},{{\bar Z}_t},{\bar \Lambda _t} + {N_{t - }})\} dt + {Z_t}d{W_t} + {{\bar Z}_t}d{M_t}\\
{P_T} =& G,
\end{split}\label{pbsre}
\end{equation}
\begin{equation}
\begin{split}
d{N_t} =&  - \{ 2({A_t}{N_{t - }} - {{\tilde \lambda }_t}{E_t}{N_{t - }}) + \left\langle {{C_t},{Z_t}} \right\rangle  + {Q_t} \\
&+ {h^ - }(t,{N_{t - }},{\Lambda _t},{{\bar \Lambda }_t},{{\bar Z}_t} + {P_{t - }})\} dt + {\Lambda _t}d{W_t} + {{\bar \Lambda }_t}d{M_t}\\
{N_T} =& G,
\end{split}\label{nbsre}
\end{equation}
where
\begin{align*}
&{h^ + }(t,p,{q_1},{q_2},{q_3}) \\
=&\mathop {\inf }\limits_{u \in \Gamma} \{ 2p{B_t}u - 2p{{\tilde \lambda }_t}{F_t}u + p{\left| {{C_t} + {D_t}u} \right|^2} + 2\left\langle {{D_t}u,{q_1}} \right\rangle  \\
&+ \left\langle {{R_t}u,u} \right\rangle  + 2({q_2} + p){{\tilde \lambda }_t}f_t^ + (1,{F_t}u)
 + 2{q_3}{{\tilde \lambda }_t}f_t^ - (1,{F_t}u)\}
\end{align*}
and
\begin{align*}
&{h^ - }(t,p,{q_1},{q_2},{q_3}) \\
=&\mathop {\inf }\limits_{u \in \Gamma} \{  - 2p{B_t}u + 2p{{\tilde \lambda }_t}{F_t}u + p{\left| { - {C_t} + {D_t}u} \right|^2} - 2\left\langle {{D_t}u,{q_1}} \right\rangle  \\
&+ \left\langle {{R_t}u,u} \right\rangle  + 2q_3 {{\tilde \lambda }_t}f_t^ + ( - 1,{F_t}u)
 + 2({q_2}+p){{\tilde \lambda }_t}f_t^ - ( - 1,{F_t}u)\}.
\end{align*}
\subsection{the Solvability of the Equations}
~~~~We have the following definitions on the solutions of the equations.
\begin{defi}We say that a pair of stochastic processes $(P,Z,\bar Z) \in \mathcal L^\infty([0,T]\times \Omega, \mathcal P(\mathbb G))\times \mathcal L^2([0,T]\times \Omega, \mathcal P(\mathbb G))\times\mathcal L^2([0,T]\times \Omega, \mathcal P(\mathbb G))$ is a solution to BSDE (\ref{pbsre}) if it satisfies the equation in the It\^o sense as well as the terminal condition and the constraint that $R+PD'D>0$. A solution $(P,Z,\bar Z)$ is called positive (resp. nonnegative) if $P>0$(resp. $P \ge 0$) and called uniformly positive if $P \ge c >0$. These definitions extent in the obvious way to the solutions of the BSDEs defined in the rest part of the paper.
\end{defi}
Before we solve the equation, let us emphasize some properties of $h^{\pm}$. First, it is obvious that
\begin{equation}
h^{\pm}(t,p,q_1,q_2,q_3) \le p |C_t|^2+ 2({q_2}+p){\tilde \lambda }_s.\label{h.estimate1}
\end{equation}
Assume that $p,q_2+p,q_3 \ge0$, we see that
\begin{equation}
\begin{split}
h^{\pm} (t,p,{q_1},{q_2},{q_3}) \ge \inf \limits_u \big\{ C\{ &\left\langle ({R_t} + pD_t'{D_t})u,u \right\rangle  - (\left| p \right| + \left| q_1 \right|)\left| u \right|\} \\
& + p{\left| {C_t} \right|^2} - (p + {q_2}){\tilde \lambda _t}\big\}
 \label{h.estimate2}
\end{split}
\end{equation}
Moreover, if $|p|,|q_1|,|q_2|,|q_3| \le n$, by (\ref{h.estimate2}), the infimum will be obtained in a bounded subset of $\Gamma$, hence is in fact a minimum  and $h^{\pm}$ are continuous with respect to $(p,q_1,q_2,q_3)$ in this situation.

Note that we get a multidimensional BSDE with quadratic growth in $z$. In general, there may be no solution for the system. See Hu and Tang \cite{HT1} for an existence result and more details on this subject. To solve the equation, we use the approach originated by Ankirchner et al \cite{ABE1} and further developed by Kharroubi and Lim \cite{KL1}: one can explicitly construct a solution by combining solutions of an associated family of Brownian BSDEs. Fortunately, we will see that we can solve these equations separately. To illustrate the idea, we give a simple example taken from \cite{KL1}. Consider the following BSDE:
\begin{equation}
\begin{split}
-dY_t&=f(U_t)dt-U_tdL_t,\qquad 0 \le t \le T,\\
Y_T&=c1_{T<\tau}+h(\tau)1_{\tau \le T}.\label{examp}
\end{split}
\end{equation}
To solve it, we first solve a recursive system of Brownian BSDEs:
\begin{align*}
Y^1_t(\theta)=h(\theta)+f(0)(T-t), \qquad \theta \wedge T \le t \le T,\\
Y^0_t=c+\int_t^Tf(Y^1_s(s)-Y^0_s)ds,\qquad 0\le t \le T.
\end{align*}
Define the process $(Y,U)$ by
\begin{align*}
Y_t=Y^0_t1_{t < \tau}+Y^1_t(\tau)1_{t \ge \tau}, \qquad 0 \le t \le T,\\
U_t=(Y^1_t(t)-Y^0_t)1_{t \le \tau}, \qquad 0 \le t \le T.
\end{align*}
By It\^o formula, we have
\begin{align*}
dY_t=&(1-L_t)dY^0_t+L_tdY_t-Y^0_tdL_t+Y^1_tdL_t\\
=&(1-L_t)f(Y^1_t(t)-Y^0_t)dt+(Y^1_t(t)-Y^0_t)dL_t\\
=&f(U_t)dt+U_tdL_t.
\end{align*}
It is also easy to see that $Y_T$ also satisfies the terminal condition. Thus $(Y,U)$ we define is a solution to (\ref{examp}).

Note that such a method is still valid in more complicate situations (see \cite{KL1} and Theorem \ref{thm1} below).We first decompose the BSDEs into two parts: the before default part and the after default part. Thus we have the following BSDEs:
\begin{equation}
\begin{split}
dP^1_t{(\theta )}= & - \{ 2A^1_t{(\theta )}P^1_t{(\theta )} + \left\langle {C^1_t{{(\theta )}},Z^1_t{{(\theta )}}} \right\rangle  + Q^1_t{(\theta )} \\
&+ {h^ + }(\theta )(t,P^1_t{(\theta )},Z^1_t{(\theta )})\} dt + Z^1_t{(\theta )}d{W_t},\quad \theta \le t \le T,\\
P^1_T{(\theta )} =& G^1(\theta ),
\end{split}\label{bdp}
\end{equation}
where
\[{h^ + }(\theta )(t,p,q) = \mathop {\inf }\limits_{u \in \Gamma} \{ 2pB^1_t{(\theta )}u +p {\left| {C^1_t{{(\theta )}} + D^1_t{{(\theta )}}u} \right|^2} + 2\left\langle {D{^1_t{(\theta )}}u,q} \right\rangle  + \left\langle {R^1_t{{(\theta )}}u,u} \right\rangle \} .\]
And
\begin{equation}
\begin{split}
dN^1_t{(\theta )} =&- \{ 2A^1_t{(\theta )}N^1_t{(\theta )} + \left\langle {C^1_t{{(\theta )}},\Lambda_t {{(\theta )}}} \right\rangle  + Q^1_t{(\theta )}\\
  &+ {h^ - }(\theta )(t,N^1_t{(\theta )},\Lambda^1_t {(\theta )})\} dt + \Lambda^1_t {(\theta )}d{W_t}, \quad \theta \le t \le T,\\
N^1_T{(\theta )} =& G^1(\theta ),
\end{split}\label{bdn}
\end{equation}
where
\[{h^ - }(\theta )(t,p,q) = \mathop {\inf }\limits_{u \in \Gamma} \{  - 2pB^1_t{(\theta )}u + p{\left| { - C^1_t{{(\theta )}} + D^1_t{{(\theta )}}u} \right|^2} - 2\left\langle {D_t{{(\theta )}}u,q} \right\rangle  + \left\langle {R^1_t{{(\theta )}}u,u} \right\rangle \} .\]
Moreover
\begin{equation}
\begin{split}
dP_t^0 =&  - \{ 2(A_t^0P_t^0 - {\lambda _t}E_t^0P_t^0) + \left\langle {C_t^0,Z_t^0} \right\rangle  + Q_t^0 \\
&+ h_0^ + (t,P_t^0,Z_t^0,P^1_t{(t)}-P_t^0,N^1_t{(t)})+\lambda_t( P^1_t{(t)}-P_t^0)\} dt\\
&+ Z_t^0d{W_t}, \qquad 0 \le t \le T,\\
P_T^0 =& {G^0},
\end{split}\label{adp}
\end{equation}
where
\begin{align*}
&h_0^ + (t,p,q,{l_1},{l_2}) \\
=& \mathop {\inf }\limits_{u \in \Gamma} \{ 2pB_t^0u - 2p{\lambda _t}F_t^0u + p{\left| {C_t^0 + D_t^0u} \right|^2} + 2\left\langle {D_t^0u,q} \right\rangle \\
 &\qquad + \left\langle {R_t^0u,u} \right\rangle  + ({l_1+p}){\lambda _t}\{ {(1 + E_t^0 + F_t^0u)^{ + ,2}} - 1\}  + {l_2}{\lambda _t}{(1 + E_t^0 + F_t^0u)^{ - ,2}}\}.
\end{align*}
And
\begin{equation}
\begin{split}
dN_t^0 =&  - \{ 2(A_t^0N_t^0 - {\lambda _t}E_t^0N_t^0) + \left\langle {C_t^0,\Lambda _t^0} \right\rangle  + Q_t^0\\
 &+ h_0^ - (t,N_t^0,\Lambda _t^0,P^1_t{(t)},N^1_t{(t)}-N_t^0)+\lambda _t(N^1_t(t)-N^0_t)\} dt \\
&+ \Lambda _t^0d{W_t},\qquad 0 \le t \le T,\\
N_T^0 =& {G^0},
\end{split}\label{adn}
\end{equation}
where
\begin{align*}
&h_0^ - (t,p,q,{l_1},{l_2})\\
 =& \mathop {\inf }\limits_{u \in \Gamma} \{  - 2pB_t^0u + 2p{\lambda _t}F_t^0u +p {\left| { - C_t^0 + D_t^0u} \right|^2} - 2\left\langle {D_t^0u,q} \right\rangle \\
 & + \left\langle {R_t^0u,u} \right\rangle  + {l_1}{\lambda _t}{(-1 - E_t^0 + F_t^0u)^{ + ,2}} + ({l_2+p}){\lambda _t}\{ {(-1 - E_t^0 + F_t^0u)^{ - ,2}} - 1\}\}.
\end{align*}
Note that we have
\[{h^ + }(t,p,{q_1},{q_2},{q_3}) = {h^+_0}(t,p,{q_1},{q_2} ,{q_3}){1_{\{ t < \tau \} }} + {h^ + }(\tau )(t,p,{q_1})1_{\tau \le t}\]
and
\[{h^ - }(t,p,{q_1},{q_2},{q_3}) = {h^-_0}(t,p,{q_1},{q_3},{q_2}){1_{\{ t < \tau \} }} + {h^ - }(\tau )(t,p,{q_1})1_{\tau \le t}.\]
We use the following theorem from \cite{KL1}
\begin{thm}\label{thm1}
Assume that for all $\theta \in \mathcal {R}_+$, the Brownian BSDEs (\ref{bdp}) and (\ref{bdn})
admit solutions $(P^1(\theta),Z^1(\theta)),(N^1(\theta),\Lambda ^1(\theta)) \in \mathcal S^\infty_{\mathbb F}[0,T]\times \mathcal L^2_{\mathbb F}[0,T]$, and that the Brownian BSDEs (\ref{adp}) and (\ref{adn}) have solutions
$(P^0,Z^0),(N^0,\Lambda ^0) \in \mathcal S^\infty_{\mathbb F}[0,T]\times \mathcal L^2_{\mathbb F}[0,T]$.
Assume moreover that $P^1(\theta)$ and $N^1(\theta)$(resp. $Z^1(\theta)$ and $\Lambda^1(\theta)$) are $\mathbb {F}\otimes \mathcal {B}(\mathcal {R}_+)$ (resp. $\mathcal {P}(\mathbb {F})\otimes \mathcal {B}(\mathcal {R}_+)$) -measurable.
If all these solutions satisfy
\[\sup \limits_{\theta}\left\|P^1(\theta)\right\|_{\mathcal S^\infty_{\mathbb F}[0,T]},\sup \limits_{\theta}\left\|N^1(\theta)\right\|_{\mathcal S^\infty_{\mathbb F}[0,T]}<\infty,\]
and
\[E[\int_{\mathcal {R}_+}\big \{\int_0^{\theta \wedge T}\left |Z^0_s\right|^2+\left |\Lambda^0_s\right|^2ds
+\int_{\theta \wedge T}^T\left |Z^1_s(\theta) \right |^2+\left |\Lambda^1_s(\theta) \right |^2ds\big \}\alpha_T(\theta)d\theta] < \infty\]
then BSDEs (\ref{pbsre}) and  (\ref{nbsre}) admit solutions $(P,Z,\bar Z),(N,\Lambda,\bar \Lambda) \in \mathcal {L}^{\infty}([0,T]\times \Omega,\mathcal {P}(\mathbb G))\times \mathcal {L}^{2}([0,T]\times \Omega,\mathcal {P}(\mathbb G))\times \mathcal {L}^{2}([0,T]\times \Omega,\mathcal {P}(\mathbb G))$ given by
\[\begin{array}{l}
\left\{ {\begin{array}{*{20}{c}}
{{P_t} = P_t^01_{t <\tau} + P^1_t(\tau)1_{\tau \le t},}\\
{{Z_t} = Z_t^01_{t \le \tau} + Z^1_t(\tau )1_{\tau < t},}\\
{{{\bar Z}_t} = (P^1_t(t ) - P_t^0)1_{t \le \tau}.}
\end{array}} \right.\\
\left\{ {\begin{array}{*{20}{c}}
{{N_t} = N_t^01_{t <\tau} + N^1_t(\tau)1_{\tau \le t},}\\
{{\Lambda_t} = \Lambda_t^01_{t \le \tau} + \Lambda^1_t(\tau )1_{\tau < t},}\\
{{{\bar \Lambda}_t} = (N^1_t(t ) - N_t^0)1_{t \le \tau}.}
\end{array}} \right.\\
\end{array}\]
\end{thm}
For the proof of this theorem, the reader can see \cite[Theorem 3.1]{KL1}.
\begin{rem}
Below, we will prove the existence of the solutions for any given $\theta$. Then we can choose $P^1$ and $N^1$ (resp. $Z^1$ and $\Lambda ^1$) as $\mathbb {F}\otimes \mathcal {B}(\mathcal {R}_+)$ (resp. $\mathcal {P}(\mathbb {F})\otimes \mathcal {B}(\mathcal {R}_+)$) -measurable processes. Indeed we know  (see Kobylanski \cite{Kob1}) that one can construct $(P^1,Z^1)$ and $(N^1,\Lambda^1)$ as limits of solutions to Lipschitz BSDEs. From Proposition C.1 in \cite{KL1},  we get $P^1$ and $N^1$ (resp. $Z^1$ and $\Lambda ^1$) as limits of $\mathbb {F}\otimes \mathcal {B}(\mathcal {R}_+)$ (resp. $\mathcal {P}(\mathbb {F})\otimes \mathcal {B}(\mathcal {R}_+)$)-measurable processes, hence also measurable.
\end{rem}
We will deal with the following two cases.\\
\begin{itshape} Standard \end{itshape}case. $Q \ge 0$,$R>0$ with $R^{-1} \in \mathcal L^{\infty}([0,T] \times \Omega,\mathcal {P}(\mathbb G),\mathcal R^{m \times m})$ and $G\ge0$.\\
\begin{itshape} Singular \end{itshape}case. $Q \ge 0$, $R \ge 0$, $G>0$ with $G^{-1} \in \mathcal L^{\infty}([0,T] \times \Omega,\mathcal {P}(\mathbb G),\mathcal R)$ and $D'D >0$ with $(D'D)^{-1} \in \mathcal L^{\infty}([0,T] \times \Omega,\mathcal {P}(\mathbb G),\mathcal R^{m \times m})$.

For the BSDE (\ref{bdp}) and (\ref{bdn}), we have the following theorem
\begin{thm}\label{thm2}
Under Assumption \ref{asm-2.2}, given any $\theta$, for the standard case, there exists a unique bounded, nonnegative maximal solution $(P^1(\theta),Z^1(\theta))$ (resp. $(N^1(\theta),\Lambda^1(\theta))$) for (\ref{bdp}) (resp. (\ref{bdn})). For the singular case, there exists a bounded, uniformly positive solution. Moreover, we have
\begin{equation}
\begin{split}
\sup \limits_{\theta}(\left \| P^1(\theta)\right \|_{\mathcal {S}^{\infty}_{\mathbb F(0,T)}}+\left \| Z^1(\theta)\right \|_{\mathcal {L}^{2}_{\mathbb F(0,T)}})< \infty,\\
\sup \limits_{\theta}(\left \| N^1(\theta)\right \|_{\mathcal {S}^{\infty}_{\mathbb F(0,T)}}+\left \| \Lambda^1(\theta)\right \|_{\mathcal {L}^{2}_{\mathbb F(0,T)}})<\infty.\label{est}
\end{split}
\end{equation}
\end{thm}
\begin{proof}
For the proof of the existence of the solutions for the extended backward Riccati equations, we refer to the paper of Hu and Zhou \cite{HZ1} (Theorem 4.1 and Theorem 4.2). Now we prove (\ref{est}).

For the standard case, we know that (see \cite{HZ1}), there exists a constant $c_1$ which only depends on the bound of the coefficents $A,B,C,D,R,G$, such that
\[\left \| P^1(\theta)\right \|_{\mathcal {S}^{\infty}_{\mathbb F(0,T)}} \le c_1.\]
Thus the norm is uniformly bounded in $\theta$. By (\ref{h.estimate2}), one can find two constants $C_1,C_2>0$ such that
\begin{equation}
p\left |C^1_t(\theta)\right |^2 \ge h^{+}(\theta)(t,p,q) \ge -\frac{C^2_1(p+\left |q\right |)^2}{C_2}+p\left |C^1_t(\theta)\right |^2.
\label{h.estimate3}
\end{equation}
Applying It\^o formula to $(P^1_t(\theta))^2$, we have
\begin{align*}
&(P^1_t(\theta))^2+\int_t^T\left|Z^1_r(\theta)\right|^2dr\\
=&(P^1_T(\theta))^2+\int^T_tP^1_r(\theta)\{2A^1_r(\theta)P^1_r(\theta)+\langle C^1_r(\theta),Z^1_r(\theta)\rangle+Q^1_r(\theta)\\
&+h^{+}(\theta)(r,P^1_r(\theta),Z^1_r(\theta))\}dr+\int^T_tP^1_r(\theta)Z^1_r(\theta)dW_r.
\end{align*}
By the boundness and non-negativity of $P$ and the inequality (\ref{h.estimate3}), taking expectation, we get that
\[\left \| Z^1(\theta)\right \|_{\mathcal {L}^{2}_{\mathbb F(0,T)}}<c_2,\]
with the constant $c_2$ independent of $\theta$. Hence we finish the proof for the standard case.

For the singular case, there will be a constant $c_3>0$ independent of $\theta$ such that
 \[c_3 \ge P^1(\theta) \ge c_3^{-1}.\]
In this case, we have
\[h^+(\theta)(t,p,q) \ge -\frac{C_1(p+\left|q\right|)^2}{pC_2}+p\left|C^1_t(\theta)\right|^2.\]
Following the same argument as above, we prove the theorem.
\end{proof}
Now we show the existence of the solution to (\ref{adp}) and (\ref{adn}). We only proof it for (\ref{adp}), since the proof is same for (\ref{adn}).
\begin{thm}\label{thm3}
Under Assumptions \ref{asm-2.1} and \ref{asm-2.2}, for the standard case there exists  a bounded, nonnegative solution $(P^0,Z^0)$ to the BSDE (\ref{adp}). And it will be uniformly positive in the singular case.
\end{thm}
\begin{proof}
For the standard case, let us first consider the following BSDE:
\begin{equation}
\begin{split}
d{P'_t} =&  - \{ 2(A_t^0 - {\lambda _t}E^0_t){P'_t} + \langle C_t^0,{Z'_t}\rangle )+(|C^0_t|^2 - {\lambda _t}){P'_t} + Q_t^0\\
&+(1+E^0_t)^2{\lambda _t}P_t^1(t)\} dt + {Z'_t}d{W_t},\\
{P'_T} =& {G^0}.
\end{split}\label{bsde-3}
\end{equation}
This is a linear BSDE with bounded coefficients and with $Q^0, P_t^1(t)\ge 0$ and $G^0 \ge 0$. Hence there exists a unique nonnegative, bounded solution $(P',Z')$. Denote by $c_1>0$ the upper bound for $P'$. Now consider the following BSDE:
\begin{equation}
\begin{split}
d{P_t} = & - \{ 2(A_t^0 - {\lambda _t}E^0_t){P_t} + \langle C_t^0,{Z_t}\rangle )- {\lambda _t}{P_t} + Q_t^0\\
 &+F(t,{P_t},{Z_t})\} dt + {Z_t}d{W_t},\\
{P_T} =& {G^0},
\end{split}\label{bsde-1}
\end{equation}
where the function $F$ is defined by
\[F(t,p,q): = [{h^0}(t,p,q,P_t^1(t) - {p},N_t^1(t))+ {\lambda _t}P_t^1(t)]{g_1}({p^ + }),\]
whereas $g_1: \mathcal {R}_+ \to [0,1]$ is a smooth truncation function satisfying $g_1(x)=1$ for $x\in[0,c_1]$, and $g_1(x)=0$ for $x\in[2c_1,+\infty)$. Note that $F$ satisfies the hypothesis (H1) of Kobylanski \cite{Kob1} thanks to the role of the truncation function $g_1$. According to \cite{Kob1}, there is a bounded maximal solution $(P,Z)$ to BSDE (\ref{bsde-1}) (see \cite{Kob1} p.565 and Theorem 2.3 for its definition and proof). Now as $F(t,p,q) \le |C^0_t|^2p+(1+E^0_t)^2{\lambda _t}P_t^1(t)$ and $(P',Z')$ is the only, hence maximal, bounded solution to (\ref{bsde-3}), we get that $P \le P' \le c_1$. Moreover, noting that $G \ge 0, Q \ge 0$ and $F(t,p,q)\ge -\frac{C_1(p+\left|q\right|)^2}{C_2}g_1(P^+)$, we conclude that $P \ge 0$ since $(0,0)$ is a solution to (\ref{bsde-1}) with $Q^0=0, G^0=0$ and $F(t,p,q)$ replaced by  $-\frac{C_1(p+\left|q\right|)^2}{C_2}g_1(P^+)$. This proves that $(P,Z)$ is a bounded nonnegative solution to (\ref{adp}).

For the singular case, we consider the following BSDE:
\begin{equation*}
\begin{split}
d{{\tilde P}_t} =&  - \{ 2(A_t^0 - {\lambda _t}E^0_t){{\tilde P}_t} + \langle C_t^0,{{\tilde Z}_t}\rangle ) - {\lambda _t}{{\tilde P}_t} + Q_t^0\\
 &+ H(t,{{\tilde P}_t},{{\tilde Z}_t})\} dt + {{\tilde Z}_t}d{W_t},\\
{{\tilde P}_T} =& {G^0},
\end{split}
\end{equation*}
where
\[H(t,p,q): =  - [p(B_t^0 - {\lambda _t}F_t^0) + (C_t^0p + q)D_t^0]{p^{ - 1}}{((D_t^0)'D_t^0)^{ - 1}}[p(B_t^0 - {\lambda _t}F_t^0) + (C_t^0p + q)D_t^0]'.\]
This is the BSDE studied in \cite{KT1} and \cite{HZ2}. By Lemma 4.1 in \cite{HZ2}, there exists a unique bounded, uniformly positive solution $(\tilde P,\tilde Z)$. Denote by $c_2$ the lower bound for $\tilde P$. Now, let us consider the following BSDE:
\begin{equation}
\begin{split}
d{P_t} = & - \{ 2(A_t^0 - {\lambda _t}E^0_t){P_t} + \langle C_t^0,{Z_t}\rangle ) - {\lambda _t}{P_t} + Q_t^0\\
 &+ \tilde F(t,{P_t},{Z_t})\} dt + {Z_t}d{W_t},\\
{P_T} =& {G^0},\label{bsde-2}
\end{split}
\end{equation}
where the function $\tilde F$ is given by
\[\tilde F(t,p,q): = [{h^0}(t,p,q,P_t^1(t) - {p},N_t^1(t)) + {\lambda _t}P_t^1(t)]{g_2}({p^ + })\]
with $g_2: \mathcal {R}_+ \to [0,1]$ being another smooth truncation function satisfying $g_2(x)=0$ for $x \in [0,\frac{1}{c_2}]$ and $g_2(x)=1$ for $x \ge c_2$. With similar discussion as in the standard case, there exists a bounded, maximal solution $(P,Z)$ of BSDE (\ref{bsde-2}).
Noting that $\tilde F(t,p,q) \ge H(t,p,q)g_2(p^+)$, the maximal solution argument gives
\[P \ge \tilde P \ge c_2.\]
This means that $(P,Z)$ is actually a bounded, uniformly positive solution to the BSDE (\ref{adp}).
\end{proof}
Combining Theorem \ref{thm1}, \ref{thm2} and \ref{thm3}, we show that there exist  bounded solutions for the system of BSDE (\ref{pbsre}) and (\ref{nbsre}).
\begin{thm}\label{thm4}
Under Assumptions \ref{asm-2.1} and \ref{asm-2.2},  either in the standard case or the singular case, there exists a bounded, nonnegative solution $(P,Z,\bar Z)$ (resp. $(N,\Lambda,\bar \Lambda)$) for the BSDE (\ref{pbsre}) (resp. (\ref{nbsre})). The solution will be uniformly positive in the singular case. Furthermore, we have that
\[{\left\| {\bar Z} \right\|_{{L^\infty }}} \le 2{\left\| P \right\|_{{L^\infty }}},{\left\| {\bar \Lambda } \right\|_{{L^\infty }}} \le 2{\left\| N \right\|_{{L^\infty }}}\]
and
\[\bar Z_t+P_{t-},\bar \Lambda _t+N_{t-} \ge 0.\]
\end{thm}

\section {Solve the constrained LQ problem}
~~~~In this section we give the optimal control for the LQ problem by the solutions to the system of BSDEs for both standard and singular case.
Define
\[\begin{array}{l}
{\xi ^ + }(t): = \mathop {\arg \min }\limits_{u \in \Gamma } {h^ + }(t,{P_{t - }},{Z_t},{{\bar Z}_t},{N_{t - }} + {{\bar \Lambda }_t}),\\
{\xi ^ - }(t): = \mathop {\arg \min }\limits_{u \in \Gamma } {h^ - }(t,{N_{t - }},{\Lambda _t},{{\bar \Lambda }_t},{P_{t - }} + {{\bar Z}_t}).
\end{array}\]
Note that the minimizers are achievable due to the discussion in the above section and $\Gamma$ is closed. By the definition, $\xi^+$ and $\xi^-$ also have the following decompositions:
\[\xi^{\pm}(t)=\xi^{\pm}_0(t)1_{t \le \tau}+\xi^{\pm}_1(t,\tau)1_{\tau < t}.\label{control.dec}\]
\begin{thm}
In both the standard and singular cases, let
$(P,Z,\bar Z),(N,\Lambda ,\bar \Lambda ) \in {\mathcal {L}^\infty }([0,T] \times \Omega ,\mathcal {P}(\mathbb G)) \times {\mathcal {L}^2}([0,T] \times \Omega ,\mathcal {P}(\mathbb G)) \times {\mathcal {L}^\infty }([0,T] \times \Omega ,\mathcal {P}(\mathbb G))$ be the bounded, nonnegative solutions to BSDEs (\ref{pbsre}) and (\ref{nbsre}) (uniformly positive in singular case). Then the following state feedback control
\begin{equation}
u^*(t) = {\xi ^ + }(t)X_t^ +  + {\xi ^ - }(t)X_t^ - \label{optim.control}
\end{equation}
is the optimal control for the LQ problem. Moreover, the value function is
\[V(t,x) = \frac{1}{2}{P_t}{x^{ + ,2}} + \frac{1}{2}{N_t}{x^{ - ,2}}.\]
\end{thm}
\begin{proof}
Now consider the state feedback control:
\begin{equation}
\begin{split}
d{X_s} =& \{ {A_s}{X_{s - }} + {B_s}({\xi ^ + }(s)X_{s - }^ +  + {\xi ^ - }(s)X_{s - }^ - )\} ds \\
&\quad+ \{ {C_s}{X_{s - }} + {D_s}({\xi ^ + }(s)X_{s - }^ +  + {\xi ^ - }(s)X_{s - }^ - )\} d{W_s}\\
 &\qquad+ \{ {E_s}{X_{s - }} + {F_s}({\xi ^ + }(s)X_{s - }^ +  + {\xi ^ - }(s)X_{s - }^ - )\} d{M_s},\\
{X_t} =& x.
\end{split}\label{optim.sde}
\end{equation}
By the lemma that follows, this equation has a c\`adl\`ag (left limit right continuous) solution. Let $(u,X)$ be any admissible control and its the corresponding state process and $(u^*,X^*)$ the state feedback control (\ref{optim.control}) and the state process. Following the discussion in Section 3, we see that the Lebesgue integrands in (\ref{value.fun}) are always positive. Define the following stopping time $\kappa _n$:
\[{\kappa _n} = \inf \{ s;\int_t^s {{{\left| {{X_r}} \right|}^2}}  + {\left| {{X_ru_r}} \right|^2} + {\left| {X_r^{ + ,2}{Z_r}} \right|^2} + {\left| {X_r^{ - ,2}{\Lambda _r}} \right|^2}dr \ge n\}  \wedge T.\]
Obviously, $\kappa_n$ is an increasing sequence of stopping time and converging to $T$ almost surely. Hence taking integration from $t$ to $\kappa _n$ and then taking conditional expectation in (\ref{value.fun}), we have
\[E^{\mathcal G_t}[\frac{1}{2}P_{\kappa_n}X^{+,2}_{\kappa_n}+\frac{1}{2}N_{\kappa_n}X^{-,2}_{\kappa_n}+\int_{t}^{\kappa _n}\{Q_rX^2_r+\langle R_ru_r,u_r\rangle\}dr]\ge\frac{1}{2}P_{t}x^{+,2}+\frac{1}{2}N_{t}x^{-,2}.\]
Letting $n \to \infty$ and noting that the processes $P$ and $N$ are quasi-left continuous, we get, from the dominated convergence theorem, that
\[J(t,x,u)\ge \frac{1}{2}P_{t}x^{+,2}+\frac{1}{2}N_{t}x^{-,2}.\]
We are now going to prove that $u^* \in {\mathcal {L}^2}([0,T] \times \Omega ,\mathcal {P}(\mathbb G))$. Once we prove this, the analysis above shows that
\[J(t,x,u^*)= \frac{1}{2}P_{t}x^{+,2}+\frac{1}{2}N_{t}x^{-,2},\]
because the Lebesgue integrand in (\ref{value.fun}) is identically zero.

In the standard case, denote by $c$ the constant such that $R \ge cI_n$. Then we have that
\[cE[{\int_t^{{\kappa _n}} {\left| {u^*(s)} \right|} ^2}ds] \le E[\frac{1}{2}{P_t}{x^{ + ,2}} + \frac{1}{2}{N_t}{x^{ + ,2}}].\]
This implies that $u^* \in {\mathcal {L}^2}([0,T] \times \Omega ,\mathcal {P}(\mathbb G))$. For the singular case, construct a sequence of stopping time as follows
\[{\theta _n} = \inf \{ s \ge t|\int_t^s {{{\left| {X_{r - }^*} \right|}^2}}  + {\left| {{C_r}X_{r - }^* + {D_r}u_r^*} \right|^2} + {\tilde \lambda _r}{\left| {{E_r}X_{r - }^* + {F_r}u_r^*} \right|^2}dr \ge n\}  \wedge T.\]
We rewrite the equation (\ref{optim.sde}) as a kind of BSDE with a random terminal time:
\begin{equation*}
\left \{
\begin{array}{l}
dX^*_s = \{ [A - B{(D'D)^{ - 1}}D'C]X^*_{s-} + B{(D'D)^{ - 1}}D'z_s\} ds + z_sd{W_s} + \bar z_sd{M_s},\\
X^*_{{\kappa _n} \wedge {\theta _n}} = X^*_{{\kappa _n} \wedge {\theta _n}},
\end{array}\right.
\end{equation*}
with $z_s = {C_s}X^*_{s-} + {D_s}{u^*_s},\bar z_{s} = {E_s}X^*_{s-} + {F_s}{u^*_s}$\\
Denote by
$$f(s): = [A - B{(D'D)^{ - 1}}D'C]X^*_{s-} + B{(D'D)^{ - 1}}D'z_s.$$
Applying It\^o formula to $(X^*_s)^2$, we get that
$$d(X^*_s)^2 = X^*_{s-}f(s) + {z_s^2} + {{\tilde \lambda }_s}{{\bar z}_s^2}ds + X^*_{s-}z_{s}d{W_s} + X^*_{s-}\bar z_{s}d{M_s}.$$
Then as in the standard estimation for the BSDE, we have
$$E[\int_t^{{\kappa _n} \wedge {\theta _n}} {{{\left| {X^*_s} \right|}^2}}  + {\left| {z_s} \right|^2} + {{\tilde \lambda }_s}{{\bar z}^2}_sds] \le \tilde cE[|X^*|^2_{{\kappa _n} \wedge {\theta _n}}] \le \frac{{\tilde c}}{c}E[\frac{1}{2}{P_t}{x^{ + ,2}} + \frac{1}{2}{N_t}{x^{ + ,2}}].$$
Appealing to Fatou's lemma, we conclude that $X^*,z \in \mathcal {L}^{2}([0,T]\times \Omega,\mathcal {P}(\mathbb {G}))$. This in turn implies that $u^* \in \mathcal {L}^{2}([0,T]\times \Omega,\mathcal {P}(\mathbb {G}))$.
\end{proof}
\begin{lem}
The equation (\ref{optim.sde}) has a c\`adl\`ag solution.
\end{lem}

\begin{proof}
Before the proof, let us illustrate the meaning of such a SDE. First, the dynamic of $X$ is governed by a Brownian SDE. Then at the random $\tau=\theta$, a jump of $X$ is induced. The size of the jump is related to $X$ and $\theta$ the time that the jump happens. After the jump, $X$ still evolves according to a Brownian SDE, but the coefficients of the SDE may be changed based on the jump time. So we can solve the SDE by decomposing it into two parts: the before default part and the after default part.  We shall rewrite the SDE (\ref{optim.sde}) into the following form
\begin{equation}
\begin{split}
d{X_t} = &\{ {{\tilde A}_t}X_{t - }^ +  + {{\hat A}_t}X_{t - }^ - \} dt + \{ {{\tilde C}_t}X_{t - }^ +  + {{\hat C}_t}X_{t - }^ - \} d{W_t} \\
&+ \{ {{\tilde E}_t}X_{t - }^ +  + {{\hat E}_t}X_{t - }^ - \} d{L_t}
\label{optim.sde2}
\end{split}
\end{equation}
where the coefficients are
\[\begin{array}{l}
{{\tilde A}_t} = {A_t} - {{\tilde \lambda }_t}{E_t} + {B_t}\xi _t^ +  - {{\tilde \lambda }_t}{F_t}\xi _t^ + ,\\
{{\hat A}_t} =  - {A_t} + {{\tilde \lambda }_t}{E_t} + {B_t}\xi _t^ -  - {{\tilde \lambda }_t}{F_t}\xi _t^ - ,\\
{{\tilde C}_t} = {C_t} + {D_t}\xi _t^ + ,{{\hat C}_t} =  - {C_t} + {D_t}\xi _t^ - ,\\
{{\tilde E}_t} = {E_t} + {F_t}\xi _t^ + ,{{\hat E}_t} =  - {E_t} + {D_t}\xi _t^ - .
\end{array}\]
Note that ${{\tilde A}_t}$ has the following form
\[{{\tilde A}_t} = \tilde A_t^0{1_{t \le \tau }} + \tilde A_t^1(\tau ){1_{\tau  < t}}\]
with some $\mathbb F$-predictable process $ \tilde A^0$ and $\mathcal P(\mathbb F) \times \mathcal B(\mathcal R_+)$-measurable process $\tilde A^1$. This is also true for the other coefficients. We will use similar notations for the decompositions.

Now consider the following SDEs:
\[\begin{array}{l}
dX_s^0 = \{ \tilde A_s^0{(X_s^0)^ + } + \hat A_s^0{(X_s^0)^ - }\} ds + \{ \tilde C_s^0{(X_s^0)^ + } + \hat C_s^0{(X_s^0)^ - }\} d{W_s},\\
X_t^0 = x
\end{array}\]
and
\begin{equation*}
\begin{split}
dX_s^1(\theta ) =& \{ \tilde A_s^1(\theta ){(X_s^1(\theta ))^ + } + \hat A_s^1(\theta ){(X_s^1(\theta ))^ - }\} ds \\
&+ \{ \tilde C_s^1(\theta ){(X_s^1(\theta ))^ + } + \hat C_s^1(\theta ){(X_s^1(\theta ))^ - }\} d{W_s},\\
X_\theta ^1(\theta ) = &X_\theta ^0 + {{\tilde E}_\theta }{(X_\theta ^0)^ + } + {{\hat E}_\theta }{(X_\theta ^0)^ - }.
\end{split}
\end{equation*}
Each SDE has a unique continuous $\mathbb F$-adapted solution (see \cite{HZ1} Lemma 5.1). Then it is obvious that the process $X_t=X^0_t1_{t < \tau}+X^1_t(\tau)1_{\tau \le t}$ is a solution to (\ref{optim.sde2}), hence a solution to (\ref{optim.sde}).
\end{proof}
\end{section}
\begin{section}{Application to Portfolio Selection}
~~~~For simplicity, we consider a financial market consisting of a bank account and one stock.  We suppose that the Brownian motion $W$ is one dimensional and $\mathbb F$ is the filtration generated by it and satisfying the usual condition.
The value of the bank count, $S_0(t)$, satisfies an ordinary differential equation:
\[\begin{array}{l}
d{S_0}(t) = r_t{S_0}(t)dt,\\
{S_0}(0) = {s_0},
\end{array}\]
where $r_t$ is deterministic. The dynamic of the risky asset is affected by other firms, the counterparties, which may default at some random time denoted by $\tau$. When the default happens, it may induce a jump in the asset price and change the dynamic of the stock. But this asset still exists and can be traded after the default of the counterparties. More precisely, let the process $L_t$ and the filtration $\mathbb G$ be what we defined in Section 2. Before the default, the stock price is governed by the following SDE:
\[dS_t^0 = S_t^0(b_t^0dt + \sigma _t^0d{W_t}),S_0^0 = {s},\]
where the coefficients are $\mathbb F$-measurable. We denote by $S^1_t(\theta),t\ge \theta$ the price of the stock after the default if the default time is at time $\theta$. At the default time $\tau$, the price has a jump
 \[S^1_\theta(\theta)=S^0_{\tau -}(1- \gamma_\theta).\]
 After the default, there is a change of regime in the coefficients depending on the default. For example, if a downward jump on the stock price is induced at default time $\tau =\theta$, the rate of the return $b^1(\theta)$ should be smaller than the rate of return $b^0$ before the default, and this gap should be increasing when the default occurs early. The stock price is still governed by an SDE for default time $\tau=\theta$:
 \[dS_t^1(\theta ) = S_t^1(\theta )(b_t^1(\theta )dt + \sigma _t^1(\theta )d{W_t}),S_\theta ^1(\theta ) = S_{\theta  - }^0(\theta )(1 - {\lambda _\theta }).\]
 Denoting by $b$ and $\sigma$ the $\mathbb G$-predictable processes $b_t=b^0_{t}1_{t \le \tau}+b^1_{t}(\tau)1_{t >\tau}$ and $\sigma_t=\sigma^0_{t}1_{t \le \tau}+\sigma^1_{t}(\tau)1_{t >\tau}$, we rewrite the price process $S$ as
\[d{S_t} = {S_{t - }}({b_t}dt + {\sigma _t}d{W_t} + {\gamma _t}d{L_t}),{S_0} = {s_0}.\]

Consider now an invest strategy that can trade continuously in this market. This is mathematically quantified by a $\mathbb G$-predictable process $\pi$ called self-financed trading strategy. It represents the money invested in the stock at time $t$. By Remark \ref{rm-2.1}, we know it has the form $\pi_t=\pi^0_{t}1_{t \le \tau}+\pi^1_t(\tau)1_{\tau <t}$. Then the wealth process $X$ is given by
\[X_t=X^0_t1_{t <\tau}+X^1_t1_{t \ge \tau},\]
where $X^0$ is the wealth process in the default-free market, governed by
\[dX_t^0 = {r_t}{X_t} + \pi _t^0((b_t^0 - {r_t})dt + \sigma _t^0d{W_t}),X_0^0 = {x_0},\]
and $X^1(\theta)$ is the wealth process after the default at time $\tau=\theta$, governed by
\[\begin{array}{l}
dX_t^1(\theta ) = {r_t}X_t^1(\theta ) + \pi _t^1(\theta )((b_t^1(\theta ) - {r_t})dt + \sigma _t^1(\theta )d{W_t}),\\
X_\theta ^1(\theta ) = X_{\theta  - }^1(\theta ) - \pi _t^1(\theta ){\gamma _\theta }.
\end{array}\]
Thus we can rewrite the wealth process as follows:
\[\begin{array}{l}
d{X_t} = {r_t}{X_{t - }} + {\pi _t}[({b_t} - {{\tilde \lambda }_t}{\gamma _t} - {r_t})dt + {\sigma _t}d{W_t} - {\gamma _t}d{M_t}],\\
{X_0} = {x_0}.
\end{array}\]
We assume that the coefficients satisfies Assumption \ref{asm-2.2} and the admissible control is the set of all square-integrable $\Gamma$-valued $\mathbb G$-predictable processes with $\Gamma= \mathcal R_{+}$ Note that we only allow $\Gamma$-valued processes, which means that the investor cannot short sell the stock.

In the mean-variance portfolio selection problem, an investor's objective is to find an admissible control
$\pi$ such that the expected terminal wealth satisfies $E[X_T]=z$, for some $z \ge x_0e^{\int_0^Tr_sds}$, while the risk measured by the variance of the terminal wealth is minimal：
\[\text{Var}(X_T):=E[(X_T-EX_T)^2]=E[X_T^2]-z^2.\]
Mathematically, it can be formulated as the following problem parameterized by z
\begin{equation}
\left\{ {\begin{array}{*{20}{c}}
{\text{Minimize }{J_{MV}}({x_0},u): = E[X_T^2] - {z^2},}\\
{\text{subject to: }E[{X_T}] = z,u( \cdot )\text{ is admissible}.}
\end{array}} \right.\label{mv-pro}
\end{equation}
The above problem is feasible if there is at least one portfolio satisfying the constraints. It is important to know when the problem is feasible for all $z \ge x_0e^{\int_0^Tr_sds}$. It means that one can select a portfolio such that its terminal wealth in average is more than the payoff in the case that one put all the money in the bank. We have the following lemma.
\begin{lem}
If we assume that $\Gamma=\mathcal {R}_+$, then the mean variance problem is feasible for all $z \ge x_0e^{\int_0^Tr_sds}$ if and only if
\begin{equation}
E[\int_0^T(b_t-r_t-\tilde \lambda_t\gamma_t)^+dt)]>0.\label{cond1}
\end{equation}
\end{lem}
\begin{proof}
We first prove the "if" part. Define
\[M:=\{(t,\omega):b_t >r_t+\tilde \lambda _t \gamma_t\}.\]
Condition (\ref{cond1}) implies that the measure of $M$ is non-zero. Consider the following  control
\[\pi_t=(b_t-r_t+\tilde \lambda _t \gamma_t)1_M.\]
It is obvious that this control is admissible. Note that for any $\beta >0$, $\beta \pi$ is also admissible. Denote by $X^{\beta}$ the corresponding wealth process. Due to the linearity of the equation, we have $X^{\beta}_t=\bar X_t+\beta \tilde X_t$, where $\bar X_t=x_0e^{\int_0^tr_sds}$ and
$\tilde X$ is the solution of the following SDE:
 \[\begin{array}{l}
d{X_t} = {r_t}{X_{t - }} + {\pi _t}[({b_t} - {{\tilde \lambda }_t}{\gamma _t} - {r_t})dt + {\sigma _t}d{W_t} - {\gamma _t}d{M_t}],\\
{X_0} = 0.
\end{array}\]
Taking expectation, we have
\begin{equation}
E[\tilde X_T]=E[\int_0^Te^{\int_t^Tr_sds}u_t(b_t-r_t-\tilde \lambda_t\gamma_t)dt].
\label{expect1}
\end{equation}
Then $E[X^{\beta}_T]=x_0e^{\int_0^Tr_sds}+\beta E[\tilde X_T]$. Due to (\ref{cond1}), $E[\tilde X_T]>0$, so we can choose $\beta$, such that $E[X^{\beta}_T]=z$.
Conversely, suppose that the problem is feasible for every $z \ge x_0e^{\int_0^Tr_sds}$. Then for some $z$, let $\pi$ be a feasible control. We can also decompose $X_t=\bar X_t+\tilde X_t$. This leads to $E\tilde X_T >0$, which implies (\ref{cond1}) by (\ref{expect1}).
\end{proof}

Finally, an optimal portfolio to (\ref{mv-pro}) is called an efficient portfolio corresponding to $z$, the corresponding $(\text{Var}(X_T),z)$ is called an efficient point. The set of all the efficient points, with $z\ge x_0e^{\int_0^Tr_sds}$, is called an efficient frontier.
The following discussion is similar to that in \cite{HZ1}, so we are not going to give all the proof. The readers can see \cite{HZ1} for details.

To handle the constraint $E[X_T]=z$, we apply Lagrange multiplier technique. Define
\[J({x_0},u;\eta ): = E[{\left| {{X_T} - \eta } \right|^2}] - {(\eta  - z)^2}.\]
We first solve the following unconstrainted problem:
\[\left\{ {\begin{array}{*{20}{c}}
{\text{Minimize }J({x_0},u;\eta )}\\
{\text{subject to: }u( \cdot )\text{ is admissible}.}
\end{array}} \right.\]
Setting $y_t=X_t-\eta e^{-\int_t^Tr_sds}$, this is exactly the singular case of constraint LQ problem we considered in Section 4. Hence we have that the optimal value $V(x)$ is
\[V(x,\eta)=P_0(x_0-\eta e^{-\int_0^Tr_sds})^{+,2}+N_0(x_0-\eta e^{-\int_0^Tr_sds})^{-,2}-(\eta  - z)^2,\]
where $P$ and $N$ is the solutions of the following BSDEs:
\[\begin{array}{l}
d{P_t} =  - \{ 2{r_t}{P_t} + {h^ + }(t,{P_t},{Z_t},{{\bar Z}_t},{{\bar \Lambda }_t} + {N_t})\} dt + {Z_t}d{W_t} + {{\bar Z}_t}d{M_t},\\
d{N_t} =  - \{ 2{r_t}{N_t} + {h^ - }(t,{N_t},{\Lambda _t},{{\bar \Lambda }_t},{{\bar Z}_t} + {P_t})\} dt + {\Lambda _t}d{W_t} + {{\bar \Lambda }_t}d{M_t},\\
{P_T} = {N_T} = 1.
\end{array}\]
Applying It\^o formula to $P_te^{-2\int _t^Tr_sds}$ and $N_te^{-2\int _t^Tr_sds}$, we have
\[1-P_0e^{-2\int _0^Tr_sds}=-E[\int_0^Te^{-2\int _t^Tr_sds}h^+(t,t,{P_t},{Z_t},{{\bar Z}_t},{{\bar \Lambda }_t} + {N_t})]\ge 0.\]
Hence $P_0e^{-2\int _0^Tr_sds}, N_0e^{-2\int _0^Tr_sds} \le 1$. In fact, the strict inequality holds with respect to $N$. If not so, we assume that $N_0e^{-2\int _0^Tr_sds}=1$, then
$$h^-(t,{N_{t-}},{\Lambda_t},{{\bar \Lambda}_t},{{\bar Z }_t} + {P_{t-}})=0,\text{\quad dt} \mathbb P \text{-a.s.}.$$
This implies that $N_t=e^{-2\int _t^Tr_sds},\Lambda_t=0$ and $\bar \Lambda_t=0$. On the other hand, we have
$$h^-(t,{N_{t-}},0,0,{{\bar Z }_t} + {P_{t-}})<0.$$
Thus we get a contradiction which implies that
\begin{equation}
N_0e^{-2\int _0^Tr_sds} < 1. \label{gtthan1}
\end{equation}
For $z=e^{\int_0^Tr_sds}$, it is obvious that the efficient portfolio is $\pi^*=0$. If $z>e^{\int_0^Tr_sds}$, applying the duality theorem, we have
\[J_{MV}^*({x_0}): = \mathop {\inf }\limits_u {J_{MV}}({x_0},u) = \mathop {\sup }\limits_{\eta  \in R} \mathop {\inf }\limits_u J({x_0},u;\eta ).\]
If $\eta < x_0e^{\int_0^Tr_sds}$, taking derivative with respect to $\eta$, we have
\begin{align*}
\frac{\partial }{{\partial \eta }}V({x_0},\eta ) =&  - 2{e^{ - \int_0^T {{r_s}ds} }}{P_0}\left( {{x_0} - \eta {e^{ - \int_0^T {{r_s}ds} }}} \right) - 2(\eta  - z)\\
 \ge&  - 2{e^{\int_0^T {{r_s}ds} }}\left( {{x_0} - \eta {e^{ - \int_0^T {{r_s}ds} }}} \right) - 2(\eta  - z) \ge 0.
\end{align*}
Thus
\[\mathop {\sup }\limits_{\eta  \in R} V({x_0},\eta ) = \mathop {\sup }\limits_{\eta  \ge {x_0}{e^{\int_0^T {{r_s}ds} }}} V({x_0},\eta ).\]
This implies that
\[J_{MV}^*({x_0}) = \frac{{{N_0}{e^{ - 2\int_0^T {{r_s}ds} }}}}{{1 - {N_0}{e^{ - 2\int_0^T {{r_s}ds} }}}}{[z - {x_0}{e^{\int_0^T {{r_s}ds} }}]^2}.\]
\textbf{Acknowledgments}\\
The author would like to thank his advisor, Prof. Shanjian Tang from Fudan University, for the helpful comments and discussions. The author would also thank the referees of this paper for helpful comments.

\end{section}
\bibliographystyle{siam}
\bibliography{myreference1}
\end{document}